\documentclass[12pt]{amsart}
\usepackage{amssymb,amsmath,amsfonts,amssymb, latexsym, verbatim, esint,amsthm, mathrsfs}
\newcommand{\be}{\begin{equation}}
\newcommand{\ee}{\end{equation}}
\newcommand{\la}{\label}
\newcommand{\ba}{\begin{array}{l}}
\newcommand{\ea}{\end{array}}
\newcommand{\Rr}{{\mathbb R}}

\newcommand{\pa}{\partial}
\newcommand{\fr}{\frac}
\newcommand{\na}{\nabla}
\newtheorem{thm}{Theorem}

\newcommand{\beg}{\begin}
\newcommand{\ov}{\overline}
\title[Local formulas for the hydrodynamic pressure]{Local formulas for the hydrodynamic pressure and applications}
\author{Peter Constantin}
\address{Department of Mathematics, Princeton University, Princeton, NJ 08544}
\email{const@math.princeton.edu}
\usepackage{color}
\usepackage[colorlinks=true, pdfstartview=FitV, linkcolor=blue, citecolor=blue,urlcolor=blue]{hyperref}

\date{\today}

\begin{document}

\begin{abstract}
We provide local formulas for the pressure of incompressible fluids. The pressure can be expressed in terms of its average and averages of squares of velocity increments in arbitrary small neighborhoods. As application, we give a brief proof of the fact that $C^{\alpha}$ velocities have $C^{2\alpha}$ (or Lipschitz) pressures. We also give some regularity criteria for 3D incompressible Navier-Stokes equations.
\end{abstract}


\maketitle

\beg{center}{\em Dedicated to the memory of Professor Mark I. Vishik.}
\end{center}

\section{Introduction} We provide local formulas for the pressure of incompressible fluids. By this we mean expressions that compute a solution of
\[
-\Delta p = \sum_{i,j=1}^3\frac{\pa^2}{\pa{x_i}\pa{x_j}}(u_iu_j),
\]
where $u$ is a divergence-free velocity, at $x\in \Omega\subset\Rr^3$, from the spherical average of
the pressure,
\[
{\ov{p}}(x,r) = \fr{1}{4\pi r^2}\int_{|x-y|=r}p(y)dS(y),
\]
and from integrals of increments $(u_i(y)-u_i(x))(u_j(y)-u_j(x))$, for $|y-x|\le r$, with arbitrary small $r$. No knowledge of the behavior of $u$ outside a small ball is needed. The main ingredient is a kind of monotonicity equation for a modified object
\[
b(x,r) = {\ov{p}}(x,r) + \fr{1}{4\pi r^2}\int_{|x-y|=r}\left(\fr{y-x}{|y-x|}\cdot u(y)\right)^2dS(y).
\]
This allows us to express the pressure as
\[
p(x) = \beta(x,r) + \pi(x,r)
\]
where  $\beta$ is just a local average of the pressure,
\[
\beta(x,r) = \fr{1}{r}\int_r^{2r}{\ov{p}}(x,\rho)d\rho,
\]
and $\pi(x,r)$ is given by a couple of integrals (\ref{pixr}) of squares of increments of velocity over a ball and over an annulus of radii $2r$. Thus, we write the pressure as a sum of two local terms, one small, and the other sufficiently well-behaved. Indeed, $\beta\in L^{\infty}(\Rr^3)$ is bounded in space (for any $r$), if $u\in L^2(\Rr^3)$ (\ref{betal2}), and $\|\na\beta\|_{L^2(\Rr^3)}$ is bounded in terms of $\|u\|_{L^4(\Rr^2)}^2$ (\ref{nabetal2}). On the other hand, $\pi$ is of the order $r^2|\na u|^2$ for small $r$. 
Well-known criteria for regularity for the 3D incompressible Navier-Stokes equations in terms of the pressure (\cite{berselli}), (\cite{sereginsverak}) do exist. If the pressure would obey the bounds that $\beta$ obeys, then regularity of solutions of the 3D Navier-Stokes equations would easily follow.  Because $\pi(x,r)\to 0$ as $r\to 0$, the suggestion that $p$ obey the same bounds as $\beta$ is not unreasonable. On the other hand, bounds on $\pi$ require some smoothness of the velocity. Higher regularity in space for velocity for weak solutions of the 3D Navier-Stokes equations was obtained in (\cite{fgt}) (see also (\cite{vasseur})).  These bounds imply that $\pi(x,r)$ is small for almost all time. For instance, $\|\pi\|_{L^3(\Rr^3)}\le C(t) r^2,\;\, t-a.e.$ (\ref{pil3}), (\ref{deltin}). The problem is that in general the time integrability of $C(t)$ is too poor to conclude regularity ($C(t)^{\fr{1}{3}}$ is time integrable, whereas $C(t)$ time integrable would be sufficient for regularity.)

The organization of this paper is as follows: In the next section we present the basic calculations which lead to the formulas for the pressure. In section 3 we give ensuing bounds for $\beta$ and $\pi$. In section 4 we give a quick proof of the bounds of higher derivatives of solutions of the 3D Navier-Stokes equations in the whole space. (The paper (\cite{fgt}) deals with spatially periodic solutions). In section 5 we give two applications: the first is a simple proof of the fact that, if $u \in C^{\alpha}$, then $p\in C^{2\alpha}$ (if $2\alpha<1$; if $2\alpha> 1$ then $p$ is Lipschitz). This result was used recently in (\cite{isett}), with a proof based on the Littlewood-Paley decomposition. A different proof (closer to ours) was obtained before, but was not published (\cite{silvestre}).
The 3D Navier-Stokes equations are regular if $u\in L^{\infty}([0,T], L^3(\Rr^3))$ (\cite{escauriazasereginsverak}), (\cite{sereginsv}). 
We give as a second application, criteria of regularity for the 3D Navier-Stokes equations in terms of $\pi$. These essentially say that if we can find $r(t)$ small such that in some sense, $\pi$ is small, and if some integral of $r(t)^{-1}$ is finite, then we have regularity.\\
Some elementary calculations needed for the formulas are presented in the Appendix.
\section{Spherical averages}
We denote
\be
{\ov{f}}(x,r) = \fr{1}{4\pi r^2}\int_{|x-y|=r}f(y)dS(y) = \fint_{|\xi| =1} f(x+r\xi)dS(\xi)
\la{ovf}
\ee
where $\fint$ denotes normalized integral. We consider solutions of
\be
-\Delta p = \na\cdot(u\cdot\na u)
\la{peq}
\ee
in $\Omega\subset \Rr^3$. We assume $\na\cdot u =0$ and smoothness of $u$.
We start by computing
\[
\ba
\pa_r {\ov{p}}(x,r)\\
= \fint_{|\xi|=1}\xi\cdot\na_x p(x+ r\xi) dS(\xi)  = \fr{1}{4\pi r}\int_{|\xi|= 1}\xi\cdot\na_{\xi}p(x+ r\xi)dS(\xi) \\
=\fr{1}{4\pi r}\int_{|\xi|<1}\Delta_{\xi} p(x + r\xi)d\xi = \fr{r}{4\pi}\int_{|\xi| <1}\Delta_x p(x+r\xi)d\xi.
\ea
\]
We use the equation (\ref{peq}). We note that, in view of the incompressibility
$\na\cdot u =0$, we have
\[
\Delta p = -\pa_i\pa_j((u_i-v_i)(u_j-v_j))
\]
for any constant vector $v$. (We use summation convention, unless explicitly stated otherwise). We have thus
\[
\ba
\pa_r {\ov{p}}(x,r) = -\fr{r}{4\pi}\int_{|\xi|<1}\pa_i\pa_j((u_i-v_i)(u_j-v_j))(x+r\xi)d\xi \\= -\fr{1}{4\pi}\int_{|\xi|<1}\pa_{\xi_i}\pa_j((u_i-v_i)(u_j-v_j))(x+ r\xi)d\xi \\ = -\fr{1}{4\pi}\int_{|\xi| =1}\xi_i(\pa_{j}((u_i-v_i)(u_j-v_j))(x+r\xi)dS(\xi)\\ = -\fr{1}{4\pi r}\int_{|\xi| =1}\xi_i\pa_{\xi_j}((u_i-v_i)(u_j-v_j))(x+ r\xi)dS(\xi).
\ea
\]
So we have 
\be
r\pa_r{\ov{p}}(x,r) = -\fint_{|\xi|=1}\xi_i\pa_{\xi_j}((u_i-v_i)(u_j-v_j)(x+ r\xi)dS(\xi).
\la{rdrpr}
\ee
\beg{lemma}\la{lem1} Let $\Omega$ be an open set in $\Rr^3$, let $x\in\Omega$. Let $r<{\mbox{dist}}(x,\partial\Omega)$, and let $u$ be a divergence-free vector field in 
$C^2(\Omega)^3$. Let $v\in \Rr^3$.  Let $p$ solve (\ref{peq}) in $\Omega$.
Then
\be
\ba
\pa_r\left\{{\ov{p}}(x,r) + \fint_{|\xi|=1}|\xi\cdot(u(x+r\xi)-v)|^2dS(\xi)\right\}\\= -\fr{1}{r}\fint_{|\xi|=1}\left[3|\xi\cdot (u(x+r\xi)-v)|^2 - |u(x+r\xi)-v|^2\right]dS(\xi).
\ea
\la{drp}
\ee
\end{lemma}
\noindent{\bf Proof.} We are going to use the identities
\be
\ba
\fint_{|\xi| =1}\xi_j\pa_{\xi_j}f(x+r\xi)dS(\xi) \\ 
= r\pa_r\left[\fint_{|\xi|=1}\xi_j^2 f(x+r\xi) dS(\xi)\right]  +\fint_{|\xi|=1} (3\xi_j^2-1)f(x+r\xi) dS(\xi)
\ea
\la{idone}
\ee
valid for each $j$, (no summation of repeated indices in the formula above), and
\be
\ba
\fint_{|\xi|=1}\left(\xi_i\pa_{\xi_j} + \xi_j\pa_{\xi_i}\right)f(x+ r\xi)dS(\xi) \\
= r\pa_r\left[\fint_{|\xi|=1}2\xi_i\xi_j f(x+r\xi)dS(\xi)\right] + \fint_{|\xi|=1}6\xi_i\xi_j f(x+r\xi)dS(\xi).
\ea
\la{idtwo}
\ee
The proofs of these identities are elementary; they are given with full detail in the Appendix.
In view of (\ref{rdrpr}), the  expression we need to average is (the negative of)
\[
\ba
\xi_1\pa_{\xi_1}(w_1^2) + \xi_2\pa_{\xi_2}(w_2^2) + \xi_3\pa_{\xi_3}(w_3^2)  +
 (\xi_1\pa_{\xi_2} + \xi_2\pa_{\xi_1})(w_1w_2)\\ + (\xi_1\pa_{\xi_3} + \xi_3\pa_{\xi_1})(w_1w_3) + (\xi_2\pa_{\xi_3} + \xi_3\pa_{\xi_2})(w_2w_3)
\ea
\]
where $w = u-v$ and the expression is evaluated at $x+r\xi$. Using (\ref{idone}), (\ref{idtwo}), we group together the terms involving $r\pa_r$, and separately the ones which do not involve $r\pa_r$,  and sum.  We obtain thus from (\ref{rdrpr})
\be
\ba
r\pa_r\ov{p}(x,r) = -r\pa_r\fint_{|\xi|=1}(\xi\cdot w)^2dS(\xi) \\
-\fint_{|\xi|=1}\left[ 3(\xi\cdot w)^2- |w|^2\right]dS(\xi),
\ea
\la{drpw}
\ee
which is the same as (\ref{drp}).

\beg{lemma} Let $x\in\Omega\subset \Rr^3$, let $0<r<{\mbox{dist}}(x,\partial \Omega)$, and let $p$ solve (\ref{peq}) with divergence-free $u\in C^{2}(\Omega)^3$. Let $v\in \Rr^3$. Then
\be
\ba
p(x) + \fr{1}{3}|u(x)-v|^2 = {\ov{p}}(x,r) + \fint_{|\xi|=1}\left |\xi\cdot(u(x+r\xi)-v)\right |^2dS(\xi)  \\ +
\int _0^r\fr{d\rho}{\rho}\fint_{|\xi|=1}\left[3\left |\xi\cdot (u(x+\rho\xi)-v)\right |^2 - |u(x+\rho \xi)-v|^2\right]dS(\xi)
\ea
\la{pv}
\ee 
\end{lemma}
{\bf Proof.} This follows immediately from (\ref{drp}) by integration $\int_0^rd\rho$, noting that
\be
{\ov{p}}(x, 0) = p(x)\la{pzero}
\ee
and
\be
\lim_{r\to 0}\fint_{|\xi|=1}\left|\xi\cdot (u(x+r\xi)-v)\right |^2dS(\xi) =\fr{1}{3}\lim_{r\to 0}\fint_{|\xi|=1}|u(x+r\xi)-v|^2dS(\xi)
\la{uz}
\ee
The formula (\ref{pv}) can be specialized by choosing $v$. Before doing this, let us introduce
\be
\sigma_{ij}(\widehat{y-x}) = 3\frac{(y_i-x_i)(y_j-x_j)}{|y-x|^2}-\delta_{ij}
\la{sigma}
\ee
where
\[
{\widehat{{y-x}}} = \fr{y-x}{|y-x|}.
\]
Note that
\[
\pa_i\pa_j\left(\fr{1}{|x-y|}\right) = \fr{\sigma_{ij}(\widehat{y-x})}{|y-x|^3}.
\]
By choosing $v=0$ in (\ref{pv}) we obtain
\be
\ba
p(x) + \fr{1}{3} |u(x)|^2 = \\
{\ov{p}}(x,r) + \fint_{|y-x|=r}|\xi\cdot u(y)|^2dS(y) +
\fr{1}{4\pi}PV\int_{B(x,r)}\fr{\sigma_{ij}(\widehat{x-y})}{|x-y|^3}(u_iu_j)(y)dy.
\ea
\la{prep}
\ee
\beg{rem}\la{rem1}
If $\Omega = \Rr^3$, if  we integrate $R^{-1}\int_R^{2R}dr$ and let $R\to\infty$ in (\ref{prep}) we obtain (assuming that $R^{-1}\int_R^{2R}\ov{p}dr$ decays)
\be
p(x) + \fr{1}{3}|u(x)|^2 = \fr{1}{4\pi}PV\int_{\Rr^3}\fr{\sigma_{ij}(\widehat{x-y})}{|x-y|^3}(u_iu_j)(y)dy
\la{puij}
\ee
a fact that follows also from
\[
p(x) = \fr{1}{4\pi}\int_{\Rr^3}\fr{1}{|x-y|}\pa_i\pa_j(u_iu_j)(y)dy
\]
by integration by parts. So (\ref{prep}) is a local version of this, valid for any $r>0$.  
\end{rem}

By choosing $v=u(x)$ in (\ref{pv}), we obtain
\be
\ba
p(x) - {\ov{p}}(x,r) - \fint_{|y-x|=r}|\xi\cdot (u(y)-u(x))|^2dS(y) \\=
\fr{1}{4\pi}\int_{B(x,r)}\fr{\sigma_{ij}(\widehat{x-y})}{|x-y|^3}((u_i(y)-u_i(x))(u_j(y)-u_j(x))dy
\ea
\la{prepx}
\ee
In order to clarify the relationship between (\ref{prep}) and (\ref{prepx}) let us observe that 
\be
\fint_{|y-x|=r}\xi_i\left(\xi\cdot u(y)\right )dS(y) + \fr{1}{4\pi}PV\int_{B(x,r)}\fr{\sigma_{ij}(\widehat{x-y})}{|x-y|^3}u_j(y)dy = \fr{1}{3}u_i(x)
\la{usig}.
\ee
This follows from the obvious fact that
\[
\fr{1}{4\pi}\int_{B(x,r)}\fr{y_i-x_i}{|y-x|^3} (\na\cdot u)(y)dy = 0
\]
by integration by parts.  
\beg{rem}\la{rem2}
Letting $r\to\infty$ we deduce from (\ref{usig}) in the whole space case, if $u$ decays, that
\be
\fr{1}{4\pi}PV\int_{\Rr^3}\fr{\sigma_{ij}(\widehat{x-y})}{|x-y|^3}u_j(y)dy = \fr{1}{3}u_i(x)
\la{sigut}
\ee
a fact that follows also from the fact that ${\mathbb P}u =u$
where ${\mathbb P}$ is the projector on divergence-free functions, using the formula
\[
{\mathbb P}v = \fr{2}{3}v +\fr{1}{4\pi}PV\int_{\Rr^3}\fr{\sigma_{ij}(\widehat{x-y})}{|x-y|^3}v_j(y)dy.
\]
\end{rem}

We write now in the principal value integral in (\ref{prep})
\[
\ba
u_i(y)u_j(y) =
(u_i(y)-u_i(x))(u_j(y)-u_j(x)) \\ + u_i(x)u_j(y) + u_j(x)u_i(y) -u_i(x)u_j(x)
\ea
\]
and take advantage of the fact that averages of $\fr{\sigma_{ij}(\widehat{y-x})}{|y-x|^3}$ on spheres centered at $x$ vanish. Using (\ref{usig}) we obtain
\[
\ba
p(x) + \fr{1}{3} |u(x)|^2 = \\
{\ov{p}}(x,r) + \fint_{|y-x|=r}|\xi\cdot u(y)|^2dS(y) +\\
\fr{1}{4\pi}PV\int_{B(x,r)}\fr{\sigma_{ij}(\widehat{x-y})}{|x-y|^3}(u_i(y)-u_i(x))(u_j(y)-u_j(x))dy\\
-2\fint_{|y-x|=r}(\xi\cdot u(x))(\xi\cdot u(y))dS(y) + \fr{2}{3}|u(x)|^2
\ea
\]
Rearranging, and noting that 
\[
\fint_{|y-x|=r}(\xi\cdot u(x))^2dS(y) = \fr{1}{3}|u(x)|^2
\]
we obtain 
\be
\ba
p(x) = \\
{\ov{p}}(x,r) + \fint_{|y-x|=r}|\xi\cdot (u(y)-u(x))|^2dS(y) +\\
\fr{1}{4\pi}\int_{B(x,r)}\fr{\sigma_{ij}(\widehat{x-y})}{|x-y|^3}(u_i(y)-u_i(x))(u_j(y)-u_j(x))dy
\ea
\la{prepu}
\ee
We have thus:
\beg{rem}\la{rem3}
The formula (\ref{prepx}) follows directly from (\ref{prep}) by using
the formula (\ref{usig}), which is a consequence of the divergence-free condition.
\end{rem}
\beg{rem}
The situation in $\Rr^2$ is entirely similar. Instead of (\ref{idone}) and (\ref{idtwo}), we have for fixed $j=1,2,$
\be
\ba
\fint_{{\mathbb S}^1}\xi_j\pa_{\xi_j} f(x+r\xi)dS(\xi) \\ =
r\pa_r\fint_{{\mathbb S}^1}\xi_j^2f(x+r\xi)dS(\xi) +\fint_{{\mathbb S}^1} (2\xi_j^2-1)f(x+r\xi) dS(\xi),
\ea
\la{id2done}
\ee
and
\be
\ba
\fint_{{\mathbb S}^1}(\xi_1\pa_{\xi_2} + \xi_2\pa_{\xi_1}) f(x+r\xi)dS(\xi) \\ = r\pa_r\fint_{{\mathbb S}^1}2\xi_1\xi_2f(x+r\xi)dS(\xi) +\fint_{{\mathbb S}^1} 2\xi_1\xi_2f(x+r\xi) dS(\xi),
\ea
\la{id2two}
\ee
and consequently, we have instead of (\ref{drpw})
\be
\ba
r\pa_r\ov{p}(x,r) = -r\pa_r\fint_{|\xi|=1}(\xi\cdot w)^2dS(\xi) \\
-\fint_{|\xi|=1}\left[ 2(\xi\cdot w)^2- |w|^2\right]dS(\xi),
\ea
\la{drp2dw}
\ee
where $w = u(x+r\xi)-v$ and $v$ is a constant vector. This again leads to a local representation formula
\be
\ba
p(x) + \fr{1}{2}|u(x)-v|^2 = {\ov{p}}(x,r) + \fint_{|\xi|=1}\left |\xi\cdot(u(x +r\xi)-v)\right |^2dS(\xi)  \\ +
\int _0^r\fr{d\rho}{\rho}\fint_{|\xi|=1}\left[2\left |\xi\cdot (u(x+\rho\xi)-v)\right |^2 - |u(x+\rho \xi)-v|^2\right]dS(\xi)
\ea
\la{p2dv}
\ee
\end{rem}

We conclude this section by mentioning similar formulae for the average of the gradient of pressure. For instance, starting from the fact that $\pa_1p$ solves the  equation
\be
-\Delta \pa_1p = \pa_i\pa_j(\pa_1(u_iu_j))
\la{pa1peq}
\ee
obtained by differentiating (\ref{peq}), we arrive at
\be
\ba
\pa_r\ov{\pa_1 p} = -\pa_r\fint_{|\xi|=1}\xi_i\xi_j\left(\pa_{x_1}(u_iu_j)(x+r\xi)\right)dS(\xi)\\ -\fr{1}{r}\fint_{|\xi|=1}(3\xi_i\xi_j-\delta_{ij})\left(\pa_{x_1}(u_iu_j)(x+r\xi)\right)dS(\xi) =\\
 -\pa_r(\fr{1}{r}\fint_{|\xi|=1}\xi_i\xi_j\left(\pa_{\xi_1}(u_iu_j)(x+r\xi)\right)dS(\xi)\\ -\fr{1}{r^2}\fint_{|\xi|=1}(3\xi_i\xi_j-\delta_{ij})\left(\pa_{\xi_1}(u_iu_j)(x+r\xi)\right)dS(\xi).
\ea
\la{gradpgradu}
\ee
We can integrate by parts in (\ref{gradpgradu}), using the relations
\be
\left\{
\ba
\fint_{|\xi|=1}\xi_1\xi_2\pa_{\xi_1}f(x+r\xi)dS(\xi) = r\pa_{r}\fint \xi_1^2\xi_2 f(x+r\xi)dS(\xi) \\ + \fint_{|\xi|=1}\left(4\xi_1^2-1\right )\xi_2 f(x+r\xi)dS(\xi),\\
\fint_{|\xi|=1}\xi_1\xi_3\pa_{\xi_1}f(x+r\xi)dS(\xi) = r\pa_{r}\fint \xi_1^2\xi_3 f(x+r\xi)dS(\xi) \\ + \fint_{|\xi|=1}\left(4\xi_1^2-1\right )\xi_3 f(x+r\xi)dS(\xi),\\
\fint_{|\xi|=1}\xi_1^2\pa_{\xi_1}f(x+r\xi)dS(\xi) =
r\pa_{r}\fint \xi_1^3f(x+r\xi)dS(\xi) \\ + \fint_{|\xi|=1}\left(4\xi_1^2-2\right )\xi_1 f(x+r\xi)dS(\xi),\\
\fint_{|\xi|=1}\xi_2^2\pa_{\xi_1}f(x+r\xi)dS(\xi) =
r\pa_{r}\fint \xi_1\xi_2^2f(x+r\xi)dS(\xi) \\ + \fint_{|\xi|=1}4\xi_1\xi_2^2 f(x+r\xi)dS(\xi),\\
\fint_{|\xi|=1}\xi_3^2\pa_{\xi_1}f(x+r\xi)dS(\xi) =
r\pa_{r}\fint \xi_1\xi_3^2f(x+r\xi)dS(\xi) \\ + \fint_{|\xi|=1}4\xi_1\xi_3^2 f(x+r\xi)dS(\xi),\\
\fint_{|\xi|=1}\xi_2\xi_3\pa_{\xi_1}f(x+r\xi)dS(\xi) =
r\pa_{r}\fint \xi_1\xi_2\xi_3f(x+r\xi)dS(\xi) \\ + \fint_{|\xi|=1}4\xi_1\xi_2\xi_3 f(x+r\xi)dS(\xi)
\ea
\la{idmany}
\right.
\ee
which can be proved in a manner similar to the proofs of (\ref{idone}), (\ref{idtwo}).  After some calculations using the relations above we arrive at
\be
\ba
\pa_r\ov{\pa_1 p} = -\left[\pa^2_r + \fr{7}{r}\pa_r + \fr{8}{r^2}\right]\fint_{|\xi| =1}\xi_1(\xi\cdot u(x+r\xi))^2dS(\xi) \\ +
\fr{2}{r}\left[\pa_r + \fr{2}{r}\right]\fint_{|\xi|=1} u_1(x+r\xi)\left(\xi\cdot u(x+r\xi)\right)dS(\xi)\\ +
\fr{1}{r}\left [\pa_r + \fr{2}{r}\right]\fint_{|\xi|=1}\xi_1 |u(x+r\xi)|^2dS(\xi).
\ea
\la{grapu}
\ee
This follows because
\be
\ov{\xi_i\xi_j\pa_{\xi_1} u_iu_j} = \left [r\pa_r +4\right]\ov{\xi_1(\xi\cdot u)^2} -2\ov{u_1(\xi\cdot u)}
\la{idthree}
\ee
and
\be
\ov{\pa_{\xi_1}|u(x+r\xi)|^2} =\left[ r\pa_r +2\right]\ov{\xi_1|u|^2} 
\la{intdxi1}
\ee

\section{Representation and bounds}
We will take $\Omega=\Rr^3$ in this section. Let us consider 
\be
b(x,r) = {\ov{p}}(x,r) + \fint_{|\xi|=1}\left |\xi\cdot u(x+r\xi)\right|^2dS(\xi)
\la{b}
\ee
The equation (\ref{drp}) with $v=0$ is
\be
\pa_r b(x,r) = r^{-1} \left[{\ov{|u|^2 -3|\xi\cdot u(y)|^2}}\right](x,r)
\la{beq}
\ee
and, integrating from $r$ to infinity, and recalling (\ref{sigma}) we obtain
\be
b(x,r) = -\fr{1}{4\pi}\int_{|x-y|\ge r}\fr{\sigma_{ij}(\widehat{x-y})}{|x-y|^3}u_i(y)u_j(y) dy
\la{bxr}
\ee
\beg{prop}\la{bineqs} Let $x\in \Rr^3$, let $r>0$, let $p$ solve (\ref{peq}) in $\Omega=\Rr^3$ with divergence-free $u\in (C^{2}(\Rr^3)\cap L^2(\Rr^3))^3$. Let $b$ be defined by (\ref{b}). Then
\be
\sup_{x\in \Rr^3}|b(x,r)| \le \fr{1}{2\pi r^3}\|u\|_{L^2}^2.
\la{bl2}
\ee 
If $u\in H^1(\Rr^2)$, then
\be
\sup_{x\in \Rr^3}|b(x,r)| \le \fr{C}{2\pi r}\|\na u\|_{L^2}^2.
\la{bh1}
\ee
where $C$ is the constant of Hardy's inequality in $\Rr^3$.
\end{prop}
\beg{rem} Obviously we do not need $C^2$ regularity for $u$, but rather enough regularity for $b$ to be defined via (\ref{b}). Of course, the representation
(\ref{bxr}) requires only $u\in L^2$.
\end{rem}
\beg{rem} The corresponding local result in an open set $\Omega$ is a bound of $b(\cdot, r)$ in $L^{\infty}(dx)$ in terms of local $L^1(dx)$ bounds for $b$ and $L^2$ (or $H^1$) bounds for $u$. This is obtained in a straightforward manner, by multiplying (\ref{beq}) by an appropriate compactly supported function of $r$ and integrating in $r$.
\end{rem} 
\noindent{\bf Proof.} The proof follows directly from the inequality
\[
\left |\sigma_{ij}(\xi)u_iu_j\right|\le 2|u|^2
\]
valid for any vector $u\in\Rr^3$ and $\xi\in {\mathbb S}^2$, and from Hardy's
inequality
\[
\int_{\Rr^3}\fr{|u(y)|^2}{|x-y|^2}dy \le C \int_{\Rr^3}|\na u(y)|^2dy.
\]
Let us define now
\be
\beta (x,r) = \fr{1}{r}\int_r^{2r}{\ov{p}}(x,\rho)d\rho
\la{beta}
\ee
\beg{prop}\la{betafty} Let $x\in \Rr^3$, let $r>0$, let $p$ solve (\ref{peq}) in  $\Omega=\Rr^3$ with divergence-free $u\in (C^{2}(\Rr^3)\cap L^2(\Rr^3))^3$. Let $\beta$ be defined by (\ref{beta}). Then
\be
\sup_{x\in \Rr^3}|\beta (x,r)| \le \fr{3}{4\pi r^3}\|u\|_{L^2}^2.
\la{betal2}
\ee
If $u\in H^1(\Rr^2)$, then
\be
\sup_{x\in \Rr^3}|\beta(x,r)| \le \fr{3C}{4\pi r}\|\na u\|_{L^2}^2.
\la{betah1}
\ee
where $C$ is the constant of Hardy's inequality in $\Rr^3$.
\end{prop}
\noindent{\bf Proof.} We note that
\[
\beta (x,r) = \fr{1}{r}\int_r^{2r}(b(x,\rho) - {\ov{(\xi\cdot u)^2}}(x,\rho))d\rho
\]
The inequalities follow in straightforward manner from 
\[
\fr{1}{r}\int_r^{2r}{\ov{(\xi\cdot u)^2}}(x,\rho)d\rho = \fr{1}{4\pi r}\int_{r\le |x-y|\le 2r}\left(\fr{x-y}{|x-y|}\cdot u(y)\right)^2\fr{dy}{|x-y|^2},
\]
Proposition \ref{bineqs} and Hardy's inequality.
\beg{rem} We introduced the average $\beta (x,r)$ of $\ov{p}(x,r)$ in order to pass from the pointwise information on $b(x,r)$ (\ref{bl2}), (\ref{bh1}), to the pointwise information on $\beta(x,r)$ (\ref{betal2}), (\ref{betah1}), without requiring other bounds than $L^2$ (or $H^1$) for $u$.  
\end{rem}
Let us consider now the weight function
\be
w(\lambda) =
\left\{
\ba
1,\quad \quad{\rm{if}}\; \; 0\le \lambda \le 1,\\
2 -\lambda \quad{\rm{if}}\;\; 1\le \lambda\le 2,\\
0 \quad\quad {\rm if}\;\; \lambda \ge 2
\ea
\right.
\la{weight}
\ee  
Let us take now the representation formula (\ref{prepx}) and average in $r$.
We obtain
\beg{thm}  Let $x\in \Rr^3$, let $r>0$, let $p$ solve (\ref{peq}) in  $\Omega=\Rr^3$ with divergence-free $u\in (C^{2}(\Rr^3)\cap L^2(\Rr^3))^3$.
Then
\be
p(x) = \beta(x,r) + \pi(x,r)
\la{betapi}
\ee
with $\beta(x,r)$ given by
\be
\beta(x,r) = \fr{1}{r}\int_r^{2r}{\ov{p}}(x,\rho)d\rho
\la{betaag}
\ee
and $\pi(x,r)$ given by
\be
\ba 
\pi(x,r) =\fr{1}{4\pi r}\int_{r\le |y-x|\le 2r}\fr{1}{|y-x|^2}\left (\fr{y-x}{|y-x|}\cdot (u(y)-u(x))\right)^2dy +\\
 \fr{1}{4\pi}\int_{|x-y|\le 2r} w\left(\fr{|y-x|}{r}\right)\fr{\sigma_{ij}(\widehat{x-y})}{|x-y|^3}(u_i(y)-u_i(x))(u_j(y)-u_j(x))dy
\ea
\la{pixr}
\ee
\end{thm}
\beg{rem}
Passing to the limit $r\to\infty$ in (\ref{betapi}) we obtain
\be
p(x) = \fr{|u(x)|^2}{3} + \fr{1}{4\pi}\int_{\Rr^3}\fr{\sigma_{ij}(\widehat{z})}{|z|^3}(u_i(x+z)-u_i(x))(u_j(x+z)-u_j(x)){dz}
\la{pdiff}
\ee
This can be obtained also from (\ref{puij}) using (\ref{sigut}).
\end{rem}
\noindent{\bf Proof.} We integrate $\fr{1}{r}\int_r^{2r}d\rho$ the representation (\ref{prepx}) written as
\be
\ba
p(x) =\\ {\ov{p}}(x,\rho) + \fint_{|y-x|=\rho}|\xi\cdot (u(y)-u(x))|^2dS(y) \\
+ \int_0^{\rho}\fr{1}{l}\fint_{|y-x|=l}\left[3(\xi\cdot (u(y)-u(x)))^2 -|u(y)-u
(x)|^2\right]dS(y)
\ea
\la{prho}
\ee
and use the fact that
\[
\fr{1}{r}\int_r^{2r}\left(\int_0^\rho f(l)dl\right)d\rho = \int_0^{2r}w\left(\fr{l}{r}\right)f(l)dl.
\]

In addition to the bounds (\ref{betal2}) and (\ref{betah1}) we also have bounds that follow from Morrey inequality
\[
\int_{\Rr^3}|u(y)|^6dy \le C \left[\int_{\Rr^3}|\na u(y)|^2dy\right]^3,
\] 
the representation
\be
p = R_iR_j(u_iu_j)
\la{prirj}
\ee
of the pressure where $R_i = \pa_i(-\Delta)^{-\fr{1}{2}}$ are Riesz transforms, and the boundedness of Riesz transforms in $L^p$ spaces.
\beg{prop}\la{blpbounds}
Let $p$ the solution of (\ref{peq}) given by (\ref{prirj}). For any $q$, $1<q<\infty$ there exist constants $C_q>0$, independent of $r>0$ so that, for any $r>0$
\be
\|{\ov{p}}(\cdot, r)\|_{L^q(\Rr^3)}\le C_q\|u\|^2_{L^{2q}(\Rr^3)}
\la{prlp}
\ee
and
\be
\|\beta(\cdot, r)\|_{L^q(\Rr^3)}\le C_q\|u\|^2_{L^{2q}(\Rr^3)}.
\la{betalp}
\ee
For any $a\in [0,2)$ there exists $C_a>0$ such that
\be
\|\beta(\cdot, r)\|_{L^3(\Rr^3)}\le C_ar^{-a}\|u\|_{L^2(\Rr^3)}^a\|\na u\|_{L^2(\Rr^3)}^{2-a}.
\la{betaint}
\ee
There exists a constant $C>0$ so that
\be
\|\na{\ov{p}}(\cdot, r)\|_{L^2}\le Cr^{-1}\|u\|_{L^4(\Rr^3)}^2
\la{naprl2}
\ee
and
\be
\|\na\beta (\cdot, r)\|_{L^2}\le Cr^{-1}\|u\|_{L^4(\Rr^3)}^2
\la{nabetal2}
\ee
\end{prop}
\noindent{\bf Proof.} The bounds (\ref{betalp}) for $\beta$ follow from the bounds (\ref{prlp}) for ${\ov{p}}$ by averaging in $r$. The bounds (\ref{prlp}) follow from (\ref{prirj}) and the boundedness of Riesz transforms in $L^p$ spaces. The bounds (\ref{betaint}) follow from (\ref{betah1}), interpolation
\[
\|\beta\|_{L^3(\Rr^3)}\le \|\beta\|_{L^{\infty}(\Rr^3)}^{\fr{a}{3}}\|\beta\|_{L^{3-a}(\Rr^3)}^{1-\fr{a}{3}},
\]
the bound (\ref{betalp}) for $q =3-a$,
\[
\|\beta(\cdot,r)\|_{L^{3-a}(\Rr^3)}\le C_a\|u\|_{L^{6-2a}(\Rr^3)}^2,
\]
and interpolation combined with the Morrey inequality
\[
\|u\|_{L^{6-2a}(\Rr^3)}\le C\|u\|_{L^2(\Rr^3)}^{\fr{a}{6-2a}}\|\na u\|_{L^2(\Rr^3)}^{\fr{6-3a}{6-2a}}.
\]
The bound (\ref{nabetal2}) follows from the bound (\ref{naprl2}) by averaging in $r$. The bound (\ref{naprl2}) follows from 
\be
\|\na{\ov{p}}(\cdot, r)\|_{L^2(\Rr^3)}\le Cr^{-1}\|{\ov{p}}(\cdot,r)\|_{L^2(\Rr^3)}
\la{naprpr}
\ee
and (\ref{prlp}) at $q=2$. The bound (\ref{naprpr}) follows from Plancherel and the observation that
\be
\widehat{{\ov{p}}}(\xi,r) = \fr{\sin(r|\xi|)}{r|\xi|}{\widehat{p}}(\xi).
\la{foupr}
\ee
Indeed,
\[
\ba
\int_{\Rr^3}e^{-ix\cdot\xi}{\ov{p}}(x,r)dx = \fint_{|\omega|=1}dS(\omega)\int_{\Rr^3}e^{-ix\cdot\xi}p(x+r\omega)dx\\
=\widehat{p}(\xi)\fint_{|\omega|=1}e^{ir\xi\cdot\omega}dS(\omega)
\ea
\]
and the last integral is computed conveniently choosing coordinates so that $\xi$ points to the North pole:
\[
\fr{1}{4\pi}\int_0^{2\pi}d\phi\int_0^{\pi}d\theta e^{ir|\xi|\cos\theta}\sin\theta d\theta = \fr{\sin(r|\xi|)}{r|\xi|}.
\]

Regarding $\pi$ we have
\beg{prop}\la{pibounds} Let $\pi(x,r)$ be defined by (\ref{pixr}).
Then
\be
|\pi(x,r)|\le C\int_{|z|\le 2r}\fr{|u(x+z)-u(x)|^2}{|z|^3}dz.
\la{pineq}
\ee
Consequently \be
\|\pi (\cdot, r)\|_{L^q}\le C_q r^2\|\na u\|^2_{L^{2q}}
\la{pinorms}
\ee
holds for all $1<q\le \infty$. In particular, at $q=3$ we have, with Morrey's inequality,
\be
\|\pi(\cdot,r)\|_{L^3}\le C r^2\|\Delta u\|^2_{L^2}.
\la{pil3}
\ee
We also have
\be
\|\pi(\cdot,r)\|_{L^q(\Rr^3)}\le C_q\|u\|_{L^{2q}(\Rr^3)}^2.
\la{pilq}
\ee
\end{prop}
\noindent{\bf{Proof.}} The inequality (\ref{pineq}) is immediate from definition. In order to prove (\ref{pinorms}) we write
\[
|u(x+z)-u(x)|^2\le |z|^2\int_0^1|\na u(x+\lambda z)|^2d\lambda
\]
and changing order of integration we have
\[
\left |\int_{\Rr^3}\phi(x)dx\int_{|z|\le 2r}\fr{|u(x+z)-u(x)|^2}{|z|^3}dz\right |
\le C r^2\|\phi\|_{L^{q'}}\|\na u\|^2_{L^{2q}}
\]
which proves (\ref{pinorms}). The bounds (\ref{pilq}) follow from (\ref{betapi}), the corresponding bounds for $p$, and (\ref{betalp}).

\section{FGT bounds in the whole space}
We take the Navier-Stokes equation 
\be
\pa_t u + u\cdot\na u -\nu\Delta u + \na p = 0,
\la{ns}
\ee
with
\be
\na\cdot u = 0,\la{divz}
\ee
multiply by $\pa_t u-\nu\Delta u$ and integrate, using incompressibility:
\[
\int_{\Rr^3}\left |\pa_t u-\nu\Delta u\right|^2dx = -\int_{\Rr^3}(u\cdot\na u)(\pa_t u -\nu\Delta u)dx.
\]
Schwartz inequality gives:
\[
\int_{\Rr^3}\left |\pa_t u-\nu\Delta u\right|^2dx\le \int_{\Rr^3}|u\cdot\na u|^2dx
\]
and so
\[
\int_{\Rr^3}\left |\pa_t u-\nu\Delta u\right|^2dx\le \|u\|_{L^{\infty}}^2\|\na u\|_{L^2}^2.
\]
The inequality
\be
\|u\|_{L^{\infty}}^2\le C \|\na u\|_{L^2}\|\Delta u\|_{L^2}
\la{ufty}
\ee
is easy to prove using Fourier transform. Thus
\[
\int_{\Rr^3}\left |\pa_t u-\nu\Delta u\right|^2dx\le C\|\Delta u\|_{L^{2}}\|\na u\|_{L^2}^3.
\]
On the other hand,
\[
\int_{\Rr^3}\left |\pa_t u-\nu\Delta u\right|^2dx= \|\pa_t u\|_{L^2}^2 + \nu^2\|\Delta u\|_{L^2}^2 + \nu\fr{d}{dt}\|\na u\|_{L^2}^2
\]
and therefore
\[
\ba
\fr{d}{dt}\|\na u\|_{L^2}^2 + \nu\|\Delta u\|^2_{L^2} + \fr{1}{\nu}\|\pa_t u\|^2_{L^2} \\
\le \fr{C}{\nu}\|\Delta u\|_{L^2}\|\na u\|_{L^2}^2\le \fr{\nu}{2}\|\Delta u\|_{L^2}^2 + \fr{C}{\nu^3}\|\na u\|_{L^2}^6 
\ea
\]
Now we denote $y(t) = \|\na u (\cdot,t)\|^2_{L^2}$, pick a constant $A>0$, divide by $(A +y)^2$ and obtain
\[
-\fr{d}{dt}\left(\fr{1}{A+y}\right) + \fr{\nu\|\Delta u\|_{L^2}^2}{(A+y)^2} + \fr{\|\pa_t u\|_{L^2}^2}{\nu(A+y)^2} \le \fr{C}{\nu^3}y.
\]
Integrating in time we obtain
\[
\int_0^T\fr{\nu\|\Delta u\|_{L^2}^2}{(A+y)^2}dt + \int_0^T \fr{\|\pa_t u\|_{L^2}^2}{\nu(A+y)^2}dt \le \fr{C}{\nu^4}\|u_0\|_{L^2}^2 + \fr{1}{A}
\]
Therefore
\be
\int_0^T\fr{\|\Delta u\|_{L^2}^2}{(A+y)^2}dt\le \fr{C}{\nu^5}\|u_0\|_{L^2}^2 + \fr{1}{\nu A} = C\nu^{-4}[D + \nu^3A^{-1}]
\la{deltau}
\ee
and
\be
\int_0^T \fr{\|\pa_t u\|_{L^2}^2}{(A+y)^2}dt \le \fr{C}{\nu^3}\|u_0\|_{L^2}^2 + \fr{\nu}{A} = C\nu^{-2}[D +\nu^{-3}A^{-1}]
\la{ut}
\ee
where we put
\[
D = \fr{\|u_0\|_{L^2}^2}{\nu}.
\]

Now
\[
\int_0^T \|\Delta u\|_{L^2}^{\fr{2}{3}}dt \le \left [\int_0^T\fr{\|\Delta u\|_{L^2}^2}{(A+y)^2}dt\right]^{\fr{1}{3}}\left[\int_0^T(A +y)dt\right]^{\fr{2}{3}}
\]
and
\[
\int_0^T \|\pa_t u\|_{L^2}^{\fr{2}{3}}dt \le\left [\int_0^T\fr{\|\pa_t u\|_{\
L^2}^2}{(A+y)^2}dt\right]^{\fr{1}{3}}\left[\int_0^T(A +y)dt\right]^{\fr{2}{3}}
\]
and therefore
\[
\int_0^T\|\Delta u\|_{L^2}^{\fr{2}{3}}dt \le C\nu^{-\fr{4}{3}}\left[D +\nu^3 A^{-1}\right]^{\fr{1}{3}}\left[ D+ AT\right]^{\fr{2}{3}}
\]
and
\[
\int_0^T \|\pa_t u\|_{L^2}^{\fr{2}{3}}dt\le C\nu^{-\fr{2}{3}}\left[D +\nu^3 A^{-1}\right]^{\fr{1}{3}}\left[D + AT \right]^{\fr{2}{3}}
\]

Now $A$ is arbitrary, but a natural explicit choice is 
\[
A^2= \nu^3T^{-1}
\]
and then we have
\be
\int_0^T\|\Delta u\|_{L^2}^{\fr{2}{3}}dt \le C\nu^{-\fr{4}{3}}\left[\fr{\|u_0\|_{L^2}^2}{\nu} + T^{\fr{1}{2}}\nu^{\fr{3}{2}}\right ]
\la{deltin}
\ee
and
\be
\int_0^T \|\pa_t u\|_{L^2}^{\fr{2}{3}}dt\le C\nu^{-\fr{2}{3}}\left[\fr{\|u_0\|\
_{L^2}^2}{\nu} + T^{\fr{1}{2}}\nu^{\fr{3}{2}}\right ].
\la{patin}
\ee
Now using the inequality (\ref{ufty}) it follows immediately that
\be
\int_0^T\|u\|_{L^{\infty}}dt \le C\nu^{-1}\left[\fr{\|u_0\|_{L^2}^2}{\nu} + T^{\fr{1}{2}}\nu^{\fr{3}{2}}\right ]^{\fr{3}{4}}\left[\fr{\|u_0\|_{L^2}^2}{\nu}\right]^{\fr{1}{4}}.
\la{intufty}
\ee
Let us consider now the other terms in (\ref{ns}). We start by computing
\[
\int_{\Rr^3}\left | u\cdot\na u + \na p\right|^2dx = 
\|u\cdot\na u\|_{L^2}^2 + \|\na p\|_{L^2}^2 + 2\int_{\Rr^3}(u\cdot\na u)\cdot(\na p) dx.
\]
Now
\[
2\int_{\Rr^3}(u\cdot\na u)\cdot(\na p) dx = -2\int p{\rm{Tr}}(\na u)^2dx =
2\int_{\Rr^3}p\Delta p dx = -2\|\na p\|_{L^2}^2.
\]
Consequently
\[
0\le \int_{\Rr^3}\left | u\cdot\na u + \na p\right|^2dx = \|u\cdot\na u\|_{L^2}^2 - \|\na p\|_{L^2}^2. 
\]
On the other hand, obviously
\[
\|u\cdot\na u\|_{L^2}\le \|u\|_{L^{\infty}}\|\na u\|_{L^2}
\]
and in view of the previous result we have
\be
\int_0^T\|u\cdot\na u\|_{L^2}^{\fr{2}{3}}dt \le C\nu^{-\fr{2}{3}}\left[\fr{\|u_0\|_{L^2}^2}{\nu} + T^{\fr{1}{2}}\nu^{\fr{3}{2}}\right ]^{\fr{1}{2}}\left[\fr{\|u_0\|_{L^2}^2}{\nu}\right]^{\fr{1}{2}}
\la{unau}
\ee 
and, because of the inequality $\|\na p\|_{L^2}\le \|u\cdot\na u\|_{L^2}$, we also have
\be
\int_0^T\|\na p\|_{L^2}^{\fr{2}{3}}dt \le  C\nu^{-\fr{2}{3}}\left[\fr{\|u_0\|_{L^2}^2}{\nu} + T^{\fr{1}{2}}\nu^{\fr{3}{2}}\right ]^{\fr{1}{2}}\left[\fr{\|u_0\|_{L^2}^2}{\nu}\right]^{\fr{1}{2}}
 \la{nap}
\ee
We have thus
\beg{thm} Let $u$ be a solution of the Navier-Stokes equation on the interval $[0,T]$. Then the quantities 
$ \|u\|_{L^{\infty}(\Rr^3)}, \|\Delta u\|_{L^2(\Rr^3)}^{\fr{2}{3}}, \|\pa_t u\|_{L^2(\Rr^3)}^{\fr{2}{3}}$,  $\|u\cdot\na u\|_{L^2(\Rr^3)}^{\fr{2}{3}},  \|\nabla p\|_{L^2(\Rr^3)}^{\fr{2}{3}}$ are almost everywhere finite on the time interval $[0,T]$, and their time integrals are bounded uniformly, with bounds (\ref{deltin}, \ref{patin}, \ref{intufty}, \ref{unau}, \ref{nap}) depending only on $T$, $\|u_0\|_{L^2(\Rr^3)}$ and $\nu$.
\end{thm}  
The proof for weak solutions follows the same pattern as the proof given above for smooth solutions, except that we mollify the advecting velocity, prove the mollification-uniform bounds and deduce the result using essentially Fatou's lemma. 
For the sake of completeness, let us mention here other estimates. Interpolating
\[
\int_0^T\|\na u\|_{L^2}^2dt<\infty
\]
and
\[
\int_0^T \|\na u\|_{L^6}^{\fr{2}{3}}dt <\infty
\]
which comes from Morrey's inequality  and (\ref{deltin}) we get
\[
\|\na u\|_{L^{3}} \le C \|\na u\|_{L^6}^{\fr{1}{2}}\|\na u\|_{L^2}^{\fr{1}{2}}
\]
which then is integrable by H\"{o}lder
\[
\int_0^T\|\na u\|_{L^3}dt<\infty.
\]
Finally, we mention that, interpolating between $L^{\infty}(dt; L^2(dx))$ and $L^2(dt; L^6(dx))$ it is easy to see that $u\in L^p(dt, L^q(dx))$ for 
$q =\fr{6p}{3p-4}$ if $p\ge 2$. For $p\in [1,2]$ interpolating between $L^{2}(dt; L^6(dx))$ and $L^{1}(dt, L^{\infty}(dx))$ we get $q=\fr{3p}{p-1}$.
 
\section{Applications}
\beg{thm}
Let $u$ solve (\ref{ns}) and (\ref{divz}) in $\Rr^3$ and assume that $u$ belongs to $L^{\infty}(dt; L^2(\Rr^3))\cap L^{2q}(dt;C^{\alpha}(\Rr^3))$ for some $q\ge 1$. Then $p\in L^q(dt; C^{2\alpha}(\Rr^3))$ if $\alpha<\fr{1}{2}$. If $\alpha =\fr{1}{2}$ then $p\in L^q(dt; LiplogLip)$ where $LiplogLip$ is the class of functions with modulus of continuity $|x-y|\log(|x-y|^{-1})$. If $\alpha>\fr{1}{2}$ then $p\in L^q(dt; Lip)$ where $Lip$ is the class of Lipschitz continuous functions.
\end{thm}
\noindent{\bf Proof.}
We start with two points $x,y$ at distance $|x-y|$ and we choose $r =8|x-y|$. The representation (\ref{prepx}) implies
\be
\left\{
\ba
\left |p(x)- {\ov{p}}(x,r)\right| \le C\|u\|_{C^{\alpha}}^2 r^{2\alpha},\\
\left |p(y)- {\ov{p}}(y,r)\right| \le C\|u\|_{C^{\alpha}}^2 r^{2\alpha},
\ea
\right.
\la{pminuspr}
\ee
so, it remains to prove that
\[
\left | {\ov{p}}(x,r) - {\ov{p}}(y,r)\right|\le Cr^{2\alpha}
\]
if $2\alpha<1$ and  $C \sim \|u\|_{C^{\alpha}}^2$. (If $2\alpha = 1$ we  obtain $r\log(r^{-1})$, and if $2\alpha>1$, $r$.) In order to do so, we use (\ref{drp}) with $v= u\left(\fr{x+y}{2}\right)$ and integrate from $r$ to infinity. We obtain
\be
\ba
{\ov{p}}(x,r) = -\fint_{|\xi|=1}\left(\xi\cdot (u(x+r\xi)-v\right)^2dS(\xi) \\ +
 \fr{1}{4\pi}\int_{|x-z|\ge r} \fr{\sigma_{ij}({\widehat{(x-z)}})}{|x-z|^3}(u_i(z)-v_i)(u_j(z)-v_j)dz
\ea
\la{pxr}
\ee
and
\be
\ba
{\ov{p}}(y,r) = -\fint_{|\xi|=1}\left(\xi\cdot (u(y +r\xi)-v\right)^2dS(\xi) \\+
\fr{1}{4\pi}\int_{|y-z|\ge r} \fr{\sigma_{ij}({\widehat{(y-z)}})}{|y-z|^3}(u_i(z)-v_i)(u_j(z)-v_j)dz
\ea
\la{pyr}
\ee
Now clearly 
\[
\left |\fint_{|\xi|=1}\left(\xi\cdot (u(x +r\xi)-v\right)^2dS(\xi)\right| \le Cr^{2\alpha}\|u\|_{C^{\alpha}}^2
\]
and
\[
\left |\fint_{|\xi|=1}\left(\xi\cdot (u(y +r\xi)-v\right)^2dS(\xi)\right| \le Cr^{2\alpha}\|u\|_{C^{\alpha}}^2,
\]
so it remains to estimate
\[
\fr{1}{4\pi}\int_{|x-z|\ge r} \fr{\sigma_{ij}({\widehat{(x-z)}})}{|x-z|^3}
w_iw_jdz - \fr{1}{4\pi}\int_{|y-z|\ge r} \fr{\sigma_{ij}({\widehat{(y-z)}})}{|y-z|^3}w_iw_jdz
\]
where $w= u(y)-v$. Now, if $|x-z|\ge r$ but  $|y-z|\le r$, then  $|x-z|\le |y-z| + |x-y| \le \fr{9}{8}r$, and so
\[
\left |\fr{1}{4\pi}\int_{|x-z|\ge r, |y-z|\le r} \fr{\sigma_{ij}({\widehat{(x-z)}})}{|x-z|^3}w_iw_jdz\right| \le C\|u\|_{C^{\alpha}}^2r^{2\alpha},
\]
and similarly, if $|y-z|\ge r$, but $|x-z|\le r$, then
\[
\left |\fr{1}{4\pi}\int_{|y-z|\ge r, |x-z|\le r} \fr{\sigma_{ij}({\widehat{(y-z)}})}{|y-z|^3}w_iw_jdz\right| \le C\|u\|_{C^{\alpha}}^2r^{2\alpha}.
\]
Finally, we are left with
\[
\fr{1}{4\pi}\int_{|x-z|\ge r, |y-z|\ge r}(K_{ij}(x-z)-K_{ij}(y-z))w_iw_jdz
\]
where
\[
K_{ij}(\zeta) = \left(3\zeta_i\zeta_j|\zeta|^{-2}- \delta_{ij}\right)|\zeta|^{-3}
\]
This is now a classical situation in singular integral theory where the smoothness of the kernel is used. We observe that
\[
\left |K_{ij}(x-z)-K_{ij}(y-z)\right| \le C|x-y|\int_0^1 |z-(y+\lambda(x-y))|^{-4}d\lambda
\]
and that $|z-(y+\lambda(x-y))|\ge \fr{7}{8} r$. Thus
\[
\ba
\left |\fr{1}{4\pi}\int_{|x-z|\ge r, |y-z|\ge r}(K_{ij}(x-z)-K_{ij}(y-z))w_iw_jdz\right| \\
\le C|x-y|\int_0^1\int_{|z-x_{\lambda}|\ge \fr{7}{8}r}|z-x_{\lambda}|^{-4}|u(z)-u\left(\fr{x+y}{2}\right)|^2dzd\lambda
\ea
\]
where $x_{\lambda} = y + \lambda(x-y)$. Now, choosing $R>0$ fixed (we could choose $R=1$, but we prefer to keep dimensionally correct quantities)
\[
\ba
|x-y|\int_0^1\int_{|z-x_{\lambda}|\ge R}|z-x_{\lambda}|^{-4}|u(z)-u\left(\fr{x+y}{2}\right)|^2dzd\lambda\\
\le C|x-y|R^{-1}\|u\|_{L^{\infty}}^2.
\ea
\]
The integral on $\fr{7r}{8}\le |z-x_{\lambda}|\le R$,
\[
|x-y|\int_0^1\int_{\fr{7r}{8}\le |z-x_{\lambda}|\le R}|z-x_{\lambda}|^{-4}|u(z)-u\left(\fr{x+y}{2}\right)|^2dzd\lambda
\]
is estimated using
\[
\left |u(z)-u\left(\fr{x+y}{2}\right)\right| \le C\|u\|_{C^{\alpha}}^2(|z-x_{\lambda}|^{2\alpha} + r^{2\alpha})
\]
The resulting bound obtained by integrating on $\fr{7}{8}r\le |z-x_{\lambda}|\le R$ is
\[
C\|u\|_{C^{\alpha}}^2|x-y|\left[\fr{1}{1-2\alpha}r^{2\alpha-1} + r^{2\alpha-1}\right]
\]
if $2\alpha<1$,
\[
C\|u\|_{C^{\alpha}}^2|x-y|\left[\log\left(\fr{8R}{r}\right ) +1-\fr{r}{R}\right]
\]
if $2\alpha =1$, and
\[
C\|u\|_{C^{\alpha}}^2|x-y|\left[\fr{R^{2\alpha -1}}{2\alpha-1} + r^{2\alpha-1}\right]
\]
if $2\alpha>1$.  This concludes the proof.

We state now some criteria for regularity. We will write $\pi(x,t, r(t))$ for $\pi$ defined according to the formula (\ref{pixr}) for a time dependent $u(x,t)$ and with a time dependent $r=r(t)$. We recall that $\pi$ is small if $u$ is regular and $r$ is small.
\begin{thm}
Let $u$ be a smooth solution of the Navier-Stokes  equation on the interval $[0,T)$. 

\noindent First criterion: Assume that there exists $U>0$, $R>0$  and $0<r(t)\le R$ such that
 \be
\int_{\{x\in \Rr^3 \; |u(x,t)|\ge U\}}|u(x,t)||\pi(x,t,r(t))|^2dx \le \fr{\nu^2}{4}\int_{\Rr^3}|u(x,t)||\na u(x,t)|^2dx
\la{cond}
\ee
holds. Assume that there exists $\gamma>4$ such that
\be
\int_0^T r(t)^{-\gamma}dt <\infty.
\la{gamma}
\ee
Then
\be
u\in L^{\infty}([0,T],  L^3(\Rr^3)).
\la{ureg}
\ee

\noindent Second criterion: Assume that there exists $r(t)$ such that $\pi = \pi(x, r(t))$ satisfies
\be
\int_0^T \|\pi\|_{L^3(\Rr^3)}^2dt <\infty
\la{pil2l3}
\ee
and that, as above, there exists $\gamma>4$ such that (\ref{gamma}) holds.
Then again (\ref{ureg}) holds.
\end{thm}
\noindent{\bf Proof.} We start with the first criterion.
We consider the evolution of the $L^3$ norm of velocity:
\[
\fr{d}{3dt}\|u\|_{L^3(\Rr^3)}^3 + \nu\int_{\Rr^3}|\na u|^2|u|dx + \int_{\Rr^3}|u|(u\cdot\na p)dx \le 0
\]
We represent $p$ using the formula (\ref{betapi}) with $r=r(t)$.
We split softly the integral involving $\pi$:
\[
\ba
\int_{\Rr^3}|u| (u\cdot\na \pi) dx = \int_{\Rr^3}\phi\left(\fr{|u|}{U}\right)|u|(u\cdot\na \pi) dx\\ +
\int_{\Rr^3}\left (1-\phi\left(\fr{|u|}{U}\right)\right )|u|(u\cdot\na\pi) dx
\ea
\]
where $\phi (q)$ is a smooth scalar function $0\le \phi(q)\le 1$, supported in $ 0\le q\le 1$.
We use the bound
\[
|\na\pi(x)| \le C\int_0^1d\lambda\int_{|z|\le 2r} \fr{dz}{|z|^2}(|\na u(x+z)| + |\na u(x))|\na u(x+\lambda z)|
\]
which follows from (\ref{pixr}) by differentiation.
It follows that
\[
\left |\int_{\Rr^3}\phi\left(\fr{|u|}{U}\right)|u|(u\cdot\na \pi) dx\right |\le CU^2r \|\na u\|^2_{L^2(\Rr^3)}.
\]
We integrate by parts in the other piece:
\[
\int_{\Rr^3}\left (1-\phi\left(\fr{|u|}{U}\right )\right )|u|(u\cdot\na\pi) dx =
-\int_{\Rr^3}\pi u\cdot\na [|u|\left (1-\phi\left(\fr{|u|}{U}\right)\right )] dx
\]
When the derivative falls on $1-\phi$ we are in the $|u|\le U$ regime and we use (\ref{pilq}) and the interpolation combined to Morrey's inequality
\[
\|u\|_{L^4(\Rr^3)}^2 \le C\|u\|_{L^3(\Rr^3)}\|\na u\|_{L^2(\Rr^3)}
\]
to deduce
\[
\ba
\left |\int_{\Rr^3}\pi |u|u\cdot\na |u|) U^{-1}\phi'\left(\fr{|u|}{U}\right) dx\right|\le 
CU\|\pi\|_{L^2(\Rr^3)}\|\na u\|_{L^2(\Rr^3)}\\
\le CU\|u\|_{L^3(\Rr^3)}\|\na u\|_{L^2(\Rr^3)}^2
\ea
\]
When the derivative falls on $|u|$ we use the condition (\ref{cond}) and the Schwartz inequality:
\[
\ba
\left |\int_{\{|u(x,t)|\ge U\}}|u\cdot\na |u| (1-\phi\left(\fr{|u|}{U}\right)\pi|dx\right|\\
\le \fr{\nu}{2}\int_{\Rr^3}|u||\na u|^2dx.
\ea
\]
As to the integral involving $\beta$, we integrate by parts, and use H\"{o}lder's inequality followed by (\ref{betaint})
\[
\ba
\left | \int_{\Rr^3} \beta u\cdot\na |u|dx \right| \le
\|\beta\|_{L^3(\Rr^3)}\|u\|_{L^3(\Rr^3)}^{\fr{1}{2}}\sqrt{\int_{\Rr^3}|u||\na u|^2dx}\\
\le \fr{1}{2\nu}\|\beta\|_{L^3(\Rr^3)}^2\|u\|_{L^3(\Rr^3)} +  {\fr{\nu}{2}}\int_{\Rr^3}|u||\na u|^2dx\\
 \le C\nu^{-1}r^{-2a}\|u\|_{L^2(\Rr^3)}^{2a}\|\na u\|_{L^2(\Rr^3)}^{4-2a}\|u\|_{L^3(\Rr^3)} + {\fr{\nu}{2}}\int_{\Rr^3}|u||\na u|^2dx 
\ea
\]
By chosing $a =\fr{\gamma}{\gamma -2}$ we have $a<2$, and using Young's inequality, we see
that
\[
r^{-2a}\|\na u\|_{L^2(\Rr^3)}^{4-2a}\le C(r^{-\gamma} + \|\na u\|^2_{L^2(\Rr^3)})
\]
is time-integrable. The upshot is that the quantity $y(t) = \|u\|_{L^3(\Rr^3)}$
obeys an ordinary differental inequality
\[
y^2\fr{dy}{dt}\le C_1(t) + C_2(t)y + C_3(t)y
\]
with $C_1(t) = CU^2r\|\na u\|_{L^2(\Rr^3)}^2$, $C_2(t) = CU\|\na u\|_{L^2(\Rr^3)}^2$ and $C_3(t) = C\nu^{-1}r^{-2a}\|\na u\|_{L^2}^{4-2a}\|u\|_{L^2(\Rr^3)}^{2a}$. The positive functions $C_1(t), C_2(t)$ and $C_3(t)$
are known to be time-integrable. The interested reader can check that the inequality above is dimensionally correct, each term has dimensions of $[L]^6[T]^{-4}$.
Then it follows that
\[
\fr{y^2}{1+y}\fr{dy}{dt}\le C_1(t)+ C_2(t) + C_3(t),
\]
(no longer dimensionally correct), and after an easy integration, it follows that $y$ is bounded a priori in time. This proves the first criterion.

For the proof of the second criterion we again represent $p = \pi(x,r) + \beta(x,r)$ with $r= r(t)$ and bound the integral involving $\pi$ using straightforward integration by parts and H\"{o}lder inequalities: 
\[
\ba
\left |\int_{\Rr^3}(u\cdot\na \pi)|u|dx\right | = \left |\int_{\Rr^3}\pi (u\cdot\na |u|)dx\right |  \\
\le \fr{\nu}{2}\int_{\Rr^3}|u||\na u|^2dx + \fr{C}{\nu}\|u\|_{L^3(\Rr^3)}\|\pi\|_{L^3(\Rr^3)}^2.
\ea
\]
We bound the contribution coming from $\beta$ the same way as we did for the first criterion. The upshot is that $y(t) = \|u\|_{L^3(\Rr^3)}$ obeys
\[
y^2\fr{dy}{dt}\le C_4(t)y + C_3(t)y
\]
with  $C_4(t) = \fr{C}{\nu}\|\pi\|_{L^3(\Rr^3)}^2$ which is time-integrable by assumption.
It follows again that $y(t)$ is bounded apriori in time.

\section{Appendix}
We prove here the identities (\ref{idone}) and (\ref{idtwo}).
We introduce polar coordinates, 
\[
\ba
\xi_1 = \rho\cos\phi \sin\theta = \rho cS,\\
\xi_2 = \rho \sin\phi\sin \theta =\rho sS,\\
\xi_3 = \rho \cos\theta = \rho C
\ea
\]
where for simplicity of notation we abbreviate $s=\sin\phi$,  $S=\sin\theta$, $c =\cos\phi$, $C = \cos\theta$.  For a function on the unit sphere $\rho=1$. But in general $f(\xi)= f(\rho cS, \rho sS, \rho C)$,  and we have
\[
\ba
f_{\theta} =\pa_{\theta} f = \rho(cCf_1+ sCf_2 -Sf_3),\\
f_{\phi} = \pa_{\phi} f = \rho(-sSf_1 + cSf_2),\\
\rho f_{\rho}=\rho \pa_{\rho} f = \rho(cSf_1 + sSf_2 + Cf_3)
\ea
\]
where $\rho\pa_{\rho}f = \xi\cdot\na_{\xi} f$ and $\na_{\xi} f = (f_1, f_2, f_3)$. We note that $\rho\pa_\rho (\fr{\xi}{|\xi|}) = 0$, for $\xi\neq 0$.
 We have
\[
\ba
Cf_\theta + S\rho f_{\rho} = \rho(cf_1+sf_2)\\
C\rho f_{\rho}-Sf_{\theta} =\rho f_3
\ea
\]
and thus
\be
\ba
\rho f_1 = c(Cf_{\theta} + S\rho f_{\rho}) -\fr{s}{S} f_{\phi},\\
\rho f_2=s(Cf_{\theta} +S\rho f_{\rho}) + \fr{c}{S} f_{\phi},\\
\rho f_3 = C\rho f_{\rho}-Sf_{\theta}
\ea
\la{naexp}
\ee
We consider now $\rho =1$ and denote for simplicity $D_{\rho}= \rho\pa_\rho$. We compute first 
\[
\fint_{|\xi|=1}\xi_1\pa_1 f(x+r\xi)dS(\xi)
\]
using of course
\[
dS(\xi) = Sd\phi d\theta.
\]
We have
\[
\ba
\xi_1\pa_{\xi_1} f = cS(c(C\pa_{\theta} +SD_{\rho}) -\fr{s}{S}\pa_{\phi})f\\
=D_{\rho}(\xi_1^2 f)  + c^2SC\pa_{\theta}f -sc\pa_{\phi} f
\ea
\] 
We used the fact that on the unit sphere $\xi = \fr{\xi}{|\xi|}$ and 
$D_{\rho}(\xi) = 0$. We multiply by $S$ and integrate, integrating by parts where possible.  In view of
\[
-c^2\fr{d}{d\theta}(S^2C) = c^2S(S^2-2C^2) = c^2S(3S^2-2)
\]
and
\[
S\fr{d}{d\phi}(sc) = 2c^2S -S,
\]
the coefficients of $f$ are obtained by adding
\[
c^2S(3S^2-2) +  2c^2S -S = S(3\xi_1^2-1),
\]
and so
\[
\ba
\fint_{|\xi|=1}\xi_1\pa_{\xi_1}f dS(\xi) =\fint_{|\xi| =1}[D_{\rho}(\xi_1^2 f) + 3\xi_1^2 f - f]dS(\xi)\\
=D_{\rho}\left[\int_{|\xi| =1}\xi_1^2 fdS(\xi)\right] + \fint_{|\xi|=1}(3\xi_1^2-1)fdS(\xi)
\ea
\]
which is the first relation in (\ref{idone}). The rest of the formulas in (\ref{idone}) are proved similarly.
Indeed, 
\[
\ba
\xi_2\pa_{\xi_2}f = sS\left[s(C\pa_{\theta}+ SD_{\rho}) +\fr{c}{S}\pa_{\theta}\right]f\\
=\left[s^2S^2D_{\rho} + s^2SC\pa_{\theta} + sc\pa_{\phi}\right]f
\ea
\] 
Upon multiplication by $S$ and integration by parts in the $\pa_{\theta}$ and
$\pa_{\phi}$ terms we obtain the coefficients of $f$
\[
\ba
-s^2\fr{d}{d\theta}(S^2C) -S \fr{d}{d\phi}(sc) = s^2S(3S^2-2) + S -2c^2S\\
=s^2S(3S^2-2) -S + 2s^2S = (3\xi_2^2-1)S
\ea
\]
and therefore
\[
\ba
\fint_{|\xi|=1}\xi_2\pa_{\xi_2}fdS(\xi)=\\
D_{\rho}\left[\fint_{|\xi|=1}\xi_2^2f dS(\xi)\right] + \fint_{|\xi|=1}(3\xi_2^2-1)fdS(\xi)
\ea
\]
like above. The third term is
\[
\xi_3\pa_3 f = C(CD_{\rho} - S\pa_{\theta})f.
\]
Multiplying by $S$ and integrating by parts the $\pa_{\theta}$ term, we compute the coefficient of $f$
\[
\fr{d}{d\theta}(CS^2) = (3C^2-1)S = (3\xi_3^2-1),
\]
and therefore we obtain the last relation of (\ref{idone})
\[
\ba
\fint_{|\xi|=1}\xi_3\pa_{\xi_3}fdS(\xi)=\\
D_{\rho}\left[\fint_{|\xi|=1}\xi_3^2f dS(\xi)\right] + \fint_{|\xi|=1}(3\xi_3^2-1)fdS(\xi).
\ea
\]
We prove now similarly the relations (\ref{idtwo}). We start with the term corresponding to the indices $(1,3)$:
\[
\ba
(\xi_1\pa_{\xi_3} + \xi_3\pa_{\xi_1})f =\\ \left[
cS(CD_{\rho}-S\pa_{\theta}) +C(c(C\pa_{\theta} + SD_{\rho}) -\fr{s}{S}\pa_{\phi})\right]f =\\
\left[2cSCD_{\rho} + (cC^2-cS^2)\pa_{\theta} -\fr{Cs}{S}\pa_{\phi}\right ]f
\ea
\]
Multiplying by $S$, integrating, and integrating by parts we obtain the coefficient of $f$ via
\[
\ba
-c\fr{d}{d\theta}(S(1-2S^2)) + C\fr{d}{d\phi}(s)\\
=-c(C-6S^2C) + Cc = 6cSCS = 6\xi_1\xi_3S
\ea
\]
and so
\[
\ba
\fint_{|\xi| =1}(\xi_1\pa_{\xi_3} + \xi_3\pa_{\xi_1})fdS(\xi) \\
=\fint_{|\xi| =1}\left[2\xi_1\xi_3 D_{\rho}f + 6\xi_1\xi_3f\right]dS(\xi) \\
=D_{\rho}\left[\fint_{|\xi|=1}2\xi_1\xi_3fdS(\xi)\right] +\fint_{|\xi|=1} 6\xi_1\xi_3 fdS(\xi)
\ea
\]
which is the $(1,3)$ relation in (\ref{idtwo}). At indices $(1,2)$ we have to compute
\[
\ba
(\xi_1\pa_2 + \xi_2\pa_1) f \\= \left [cS(sSD_{\rho} +sC\pa_{\theta} +\fr{c}{S}\pa_{\phi}) + sS(cSD_{\rho} + cC\pa_{\theta} - \fr{s}{S}\pa_{\phi})\right]f\\
= 2cSsSD_{\rho}f + 2cs(SC)\pa_{\theta}f + (c^2-s^2)\pa_{\phi}f.
\ea
\]
Multiplying by $S$ and integrating by parts, we obtain the coefficient of $f$ via
\[
\ba
-2cs\fr{d}{d\theta}(S^2C) -S\fr{d}{d\phi}(c^2-s^2) =\\
2cs(S^3-2SC^2) + 4 Scs = 2cs(S^3-2S+ 2S^3) +4 csS = 6csS^3\\
= 6\xi_1\xi_2S.
\ea
\]
We obtained thus
\[
\ba
\fint_{|\xi| =1}(\xi_1\pa_{\xi_2} + \xi_2\pa_{\xi_3})fdS(\xi) \\
=\fint_{|\xi| =1}\left[2\xi_1\xi_2 D_{\rho}f + 6\xi_1\xi_2f\right]dS(\xi) \\
=D_{\rho}\left[\fint_{|\xi|=1}2\xi_1\xi_2fdS(\xi)\right] +\fint_{|\xi|=1} 6\xi_1\xi_2 fdS(\xi)
\ea
\]
which is the $(1,2)$ relation of (\ref{idtwo}). Finally, at $(2,3)$ we have to compute
\[
\ba
(\xi_2\pa_3 + \xi_3\pa_2)f = sS(CD_{\rho}-S\pa_{\theta})f + C(sSD_{\rho} +sC\pa_{\theta} + \fr{c}{S}\pa_{\phi})f \\ = 2sSCD_{\rho}f  + (s(C^2-S^2)\pa_{\theta}+ C\fr{c}{S}\pa_{\phi})f. 
\ea
\]
Multiplying by $S$ and integrating by parts, the coefficient of $f$ is computed via
\[
\ba
-s\fr{d}{d\theta}(S(C^2-S^2)) - C\fr{d}{d\phi}c =\\
s(6S^2C-C) + Cs = 6sSCS = 6\xi_2\xi_3S\\
\ea
\]
and we obtain thus
\[
\ba
\fint_{|\xi| =1}(\xi_2\pa_{\xi_3} + \xi_3\pa_{\xi_2})fdS(\xi) \\
=\fint_{|\xi| =1}\left[2\xi_2\xi_3 D_{\rho}f + 6\xi_2\xi_3f\right]dS(\xi) \\
=D_{\rho}\left[\fint_{|\xi|=1}2\xi_2\xi_3fdS(\xi)\right] +\fint_{|\xi|=1} 6\xi_2\xi_3 fdS(\xi)
\ea
\]
which is the $(2,3)$ relation of (\ref{idtwo}).

\noindent{\bf Acknowledgment} Research partially supported by grants NSF-DMS 1209394 and NSF-DMS 1265132.

\end{document}